\magnification1200
\input amssym.def 
\input amssym.tex 
\def\SetAuthorHead#1{}\def\SetTitleHead#1{}
\def\NoindentAfter{\everypar={\setbox0=\lastbox\everypar={}}}
\def\H#1\par#2\par{{\baselineskip=15pt\parindent=0pt\parskip=0pt
 \leftskip= 0pt plus.2\hsize\rightskip=0pt plus.2\hsize
 \bf#1\unskip\break\vskip 4pt\rm#2\unskip\break\hrule
 \vskip40pt plus4pt minus4pt}\NoindentAfter}
\def\HH#1\par{{\bigbreak\noindent\bf#1\medskip}\NoindentAfter}
\def\HHH#1\par{{\bigbreak\noindent\bf#1\unskip.\kern.4em}}
\def\th#1\par{\medbreak\noindent{\bf#1\unskip.\kern.4em}\it}
\def\endth{\medbreak\rm}
\def\pf#1\par{\medbreak\noindent{\it#1\unskip.\kern.4em}}
\def\df#1\par{\medbreak\noindent{\it#1\unskip.\kern.4em}}
\def\enddf{\medbreak}
\let\rk\df\let\endrk\enddf
\let\Roster\bgroup\let\endRoster\egroup
\def\\{}\def\text#1{\hbox{\rm #1}}
\def\mop#1{\mathop{\rm\vphantom{x}#1}\nolimits}
\def\MaxReferenceTag#1{}
\def\qedbox{\vrule width2mm height2mm\hglue1mm\relax}
\def\qed{\ifmmode\qedbox\else\hglue5mm\unskip\hfill\qedbox\medbreak\fi\rm}
\let\SOverline\overline

\let\Item\item\let\ItemItem\itemitem
\def\cite#1{{\bf[#1]}}
\def\Em#1{{\it #1\/}}
\def\Bib#1\par{\bigbreak\bgroup\centerline{\bf #1}\medbreak\parindent30pt
 \parskip2pt\frenchspacing\par}
\def\endBib{\par\egroup}
\newdimen\Overhang
\def\rf#1{\par\noindent\hangafter1\hangindent=\parindent
     \setbox0=\hbox{[#1]}\Overhang\wd0\advance\Overhang.4em\relax
     \ifdim\Overhang>\hangindent\else\Overhang\hangindent\fi
     \hbox to \Overhang{\box0\hss}\ignorespaces}
\def\bbH{{\Bbb H}}

\def\bbZ{{\Bbb Z}}
\def\Coordinates{\bigbreak\bgroup\parindent=0pt\obeylines}
\def\endCoordinates{\egroup}

\def\Title{Regular Cocycles and Biautomatic Structures}
\def\Authors{Walter D. Neumann and Lawrence Reeves}
\def\eval#1{\SOverline{#1}}
\def\len{\mop{len}}
\SetAuthorHead{\Authors}
\SetTitleHead{\Title}

\H \Title\par\Authors\par

In \cite{ECHLPT} it is shown that if the fundamental group of a
Seifert fibred 3-manifold is not virtually nilpotent then it has an
automatic structure.  In the unpublished 1992 preprint \cite{G2}
Gersten constructs a biautomatic structure on the fundamental group of
any circle bundle over a hyperbolic surface.  He asks if the same can
be done for the above Seifert fibered 3-manifold.  We show the
existence of such a biautomatic structure.

We do this in the context of a general discussion of biautomatic
structures on virtually central extensions of finitely generated
groups. A \Em{virtually central extension} is an extension of a group
$G$ by an abelian group $A$ for which the induced action of $G$ on $A$
is \Em{finite}, that is, given by a map $G\to\mop{Aut}(A)$ with finite
image.  The fundamental group of a Seifert fibered 3-manifold as above
is a virtually central extension of a Fuchsian group $G$ by
$\bbZ$. (For convenience we are using the term ``Fuchsian group'' for
any discrete finitely generated subgroup of $\mop{Isom}(\bbH^2)$ ---
orientable or not.)

We use a concept of ``regular 2-cocycles'' on a group $G$ which was
suggested by Gersten's work.  Here ``regularity'' is with respect to a
(possibly asynchronously) automatic structure $L$ on $G$.  If $L$ is a
biautomatic structure on $G$ we show that any virtually central
extension of $G$ defined by an $L$-regular cocycle also has a
biautomatic structure.

As an application we show that any virtually central extension of a
Fuchsian group $G$ by a finitely generated abelian group $A$ is
biautomatic.  In fact, if $L$ is a geodesic language on $G$, we show
that all of $H^2(G;A)$ is represented by $L$-regular cocycles.  In
case $G$ is torsion free and $A=\bbZ$ with trivial $G$-action this is
implicit in Gersten's work (loc.~cit.~--- we give an independent
treatment here that is more geometric; alternatively, it follows from
his result about biautomaticity plus Theorem A below).  The general
case follows easily from this using Corollary 2.7 below, which says
that a cohomology class for a group $G$ is regular if its restriction
to some finite index subgroup of $G$ is regular.

The converse to the fact that regular cocycles lead to biautomatic
structures is also true.

\th Theorem A

Let $E$ be a virtually central extension of the group $G$ by a
finitely generated abelian group $A$.  Then $E$ carries a biautomatic
structure if and only if $G$ has a biautomatic structure $L$ for which
the cohomology class of the extension is represented by an $L$-regular
cocycle.
\endth

This strengthens the result of Lee Mosher \cite{M} that biautomaticity
of a central extension of $G$ implies biautomaticity of $G$. We use
his work in the proof of Theorem A.

\HH 1.~~Basic Definitions

Let $G$ be a finitely generated group and $X$ a finite set which maps
to a monoid generating set of $G$. The map of $X$ to $G$ can be
extended in the obvious way to give a monoid homomorphism of $X^*$
onto $G$ which will be denoted by $w\mapsto\eval w$.  For convenience
of exposition we will always assume our generating sets are
\Em{symmetric}, that is, they satisfy $\eval X=\eval X^{-1}$.  If
$L\subset X^*$ then the pair consisting of $L$ and the evaluation map
$L\to G$ will be called a \Em{language on $G$}.  Abusing terminology,
we will often suppress the evaluation map and just call $L$ the
language on $G$ (but therefore, we may use two letters, say $L$ and
$L'$, to represent the same language $L\subset X^*$ with two different
evaluation maps to two different groups).  A language on $G$ is a
\Em{normal form} if it surjects to $G$.

A \Em{rational structure} for $G$ is a normal form $L\subset X^*$ for
$G$ which is a regular language (i.e., the set of accepted words for
some finite state automaton).

The \Em{Cayley graph} $\Gamma_X(G)$ is the directed graph with vertex
set $G$ and a directed edge from $g$ to $g\eval x$ for each $g\in G$
and $x\in X$; we give this edge a label $x$. 

Each word $w \in X^*$ defines a path $[0,\infty)\to\Gamma$ in the
Cayley graph $\Gamma=\Gamma_X(G)$ as follows (we denote this path also
by $w$): $w(t)$ is the value of the $t$-th initial segment of $w$ for
$t=0,\ldots,\len(w)$, is on the edge from $w(s)$ to $w(s+1)$ for
$s<t<s+1\le\len(w)$ and equals $\eval w$ for $t\ge\len(w)$. We refer
to the translate by $g\in G$ of a path $w$ by $gw$.

Let $\delta\in \Bbb N$. Two words $v,w\in X^*$ \Em{synchronously
$\delta$-fellow-travel} if the distance $d(w(t),v(t))$ never exceeds
$\delta$.  They \Em{asynchronously $\delta$-fellow-travel} if there
exist non-decreasing proper functions $t\mapsto t', t\mapsto t''\colon
[0,\infty)\to[0,\infty)$ 
such that $d({v(t')},{w(t'')})\le\delta$ for
all $t$.

A rational structure $L$ for $G$ is a \Em{synchronous} resp.\
\Em{asynchronous automatic structure} if there is a constant $\delta$
such that any two words $u,v\in L$ with $d(\eval u,\eval v)\le 1$
synchronously resp.\ asynchronously fellow-travel. A synchronous
automatic structure $L$ is \Em{synchronously biautomatic} if there is
a constant $\delta$ such that if $v,w\in L$ satisfy $\eval
w=\eval{xv}$ with $x\in X$ then $\eval xv$ and $w$ synchronously
$\delta$-fellow-travel.  See \cite{NS1} for a discussion of the
relationship of these definitions with those of \cite{ECHLPT}.  In
particular, as discussed there, if $L\to G$ is finite-to-one, then the
definitions are equivalent; by going to a sublanguage of $L$ this can
always be achieved.

We define two rational structures $L$ and $L'$ on $G$ to be
\Em{equivalent}, written $L\sim L'$, if there exists a $\delta$ such
that every $L$-word is asynchronously $\delta$-fellow-travelled by an
$L'$ word with the same value and vice versa.  If $L$ and $L'$ are
asynchronous automatic structures this is equivalent to requiring that
$L\cup L'$ be an asynchronous automatic structure.

If $L$ is a rational structure on $G$ we say a subset $S\subset G$ is
\Em{$L$-rational} if the language $$L_S := \{w\in L : \eval w \in
S\}$$ is a regular language. The subset $S$ is \Em{$L$-quasiconvex} if
there exists a $\delta$ such that every $w\in L$ with $\eval w\in S$
travels in a $\delta$-neighborhood of $S\subset \Gamma_X(G)$.  The
following is well-known (e.g., \cite{GS}, \cite{NS1}).

\th Proposition 1.1 

If $L\sim L'$ are equivalent rational structures on $G$ then any
subset $H$ of $G$ is $L$-rational if and only if it is $L'$-rational.
Moreover, if $H$ is a subgroup, then it is $L$-rational if and only if
it is $L$-quasiconvex.\qed
\endth

We shall also need the following.

\th Lemma 1.2

Let $L$ be an asynchronous automatic structure on $G$.  If $S$ is an
$L$-rational subset of $G$ then so is $Sg$ for any $g\in G$.
Moreover, if $L$ is a biautomatic structure then $gS$ is also
$L$-rational.  In particular, if $H$ is a subgroup of finite index
then its right-cosets $Hg$ are $L$-rational in the automatic case and
two-sided cosets $g_1Hg_2$ are $L$-rational in the biautomatic
case.\endth

\pf Proof

Suppose $L$ is an automatic structure and $S$ is $L$-rational.  It
suffices to show that $S\eval x$ is rational for any generator $x$.
We can use a standard comparator automaton (cf.~\cite{ECHLPT}) to see
that $\{(u,v)\in L^2:\eval v = \eval{ux},\eval u\in S\}$ is the
language of an asynchronous two-tape automaton.  The projection onto
the second factor is therefore a regular language, but it is just the
language of words $v\in L$ that evaluate into $S\eval x$.  Thus
$S\eval x$ is $L$-rational.  The proof that $gS$ is rational if $L$ is
biautomatic is completely analogous. The final sentence of the lemma
then follows since a subgroup of finite index, being quasiconvex, is
rational by Proposition 1.1.\qed

\HH 2.~~$L$-Regular Cocycles and Biautomatic Structures.

Let $G$ be a group and $A$ be a finitely generated abelian group.
Suppose
$$0\to A\mathop{\to}\limits^\iota E\mathop{\to}\limits^\pi G\to 1$$ 
is a virtually central extension of $G$.  We write $A$ additively and
we denote the action of an element $g\in G$ on $A$ by $a\mapsto
a^g$. Choose a section $s:G\to E$.  Then a general element of $E$ has
the form $s(g)\iota(a)$ with $g\in G$ and $a\in A$ and the group
structure in $E$ is given by a formula
$$s(g_1)\iota(a_1)s(g_2)\iota(a_2)
  =s(g_1g_2)\iota(a_1+a_2+\sigma(g_1,g_2)),$$
where $\sigma\colon G\times G\to A$ is a 2-cocycle on $G$ with
coefficients in the $G$-module $A$.  Changing the choice of section
changes the cocycle $\sigma$ by a coboundary.  Conversely, given a
cocycle $\sigma$, the above multiplication rule defines a virtually
central extension of $G$ by $A$.

\df Definition

Suppose $G$ has finite generating set $X$ and $L\subset X^*$ is an
asynchronous automatic structure on $G$.  We say a 2-cocycle $\sigma$
as above is \Em{weakly bounded} if 
 \Roster
 \ItemItem{1.} The sets $\sigma(X,G)$ and $\sigma(G,X)$ are finite;
\endRoster
\noindent and is \Em{$L$-regular} if in addition
 \Roster
 \ItemItem{2.} For each $x\in X$ and $a\in A$ the
subset $\{g\in G: \sigma(g,x)=a\}$ is an $L$-rational subset of $G$.
\endRoster
\noindent A cohomology class in $H^2(G;A)$ is \Em{$L$-regular} if it
can be represented by an $L$-regular cocycle.  The term ``weakly
bounded'' reflects the standard terminology of ``bounded'' for a
cocycle that satisfies $\sigma(G,G)$ finite.
\enddf

\th Lemma 2.1 

\Item{1.} If $\sigma$ is an $L$-regular cocycle then for any $h\in G$
and $a\in A$ the set $\{g\in G:\sigma(g,h)=a\}$ is an $L$-rational
subset of $G$
\Item{2.} If $L_1$ and $L_2$ are equivalent asynchronous 
automatic structures then
any $L_1$-regular cocycle is $L_2$-regular.\endth

\pf Proof

It is enough to show that if the statement of Lemma 2.1.1 is true for
$h_1$ and $h_2$ then it is true for $h=h_1h_2$.  Now the cocycle
relation says
$$\sigma(g,h_1h_2)=\sigma(g,h_1)^{h_2} +
\sigma(gh_1,h_2)-\sigma(h_1,h_2).$$  Thus
$$\eqalign{\{g\in G:\sigma(g,h_1h_2)=a\} =&\cr
\bigcup_{a_1+a_2=a+\sigma(h_1,h_2)} &\{g\in G:
\sigma(g,h_1)=a_1^{h_1^{-1}}\}\cap\{g\in G:\sigma(gh_1,h_2)=a_2\}.\cr}
$$
This is a finite union, since the sets on the right are empty for all
but finitely many values of $a_1$ and $a_2$.  It is a union of
rational subsets since $\{g\in G:\sigma(gh_1,h_2)=a_2\}=\{g\in
G:\sigma(g,h_2)=a_2\}h_1^{-1}$ is a right-translate of an $L$-rational
subset and hence $L$-rational by Lemma 1.2.  This proves 2.1.1.

Part 2 of the lemma follows from the fact that a subset of $G$ is
$L_1$-rational if and only if it is $L_2$-rational (Proposition 2.1).
Note that we also use part 1 of the lemma, since $L_1$ and $L_2$ may
be languages on different generating sets.\qed

Suppose now that $L$ is a finite-to-one biautomatic structure on $G$
and that $E$ is a virtually central extension as above given by a regular
cocycle $\sigma$ determined by a section $s$. Consider the finite
subset $\{s(x)\iota(-\sigma(g,x)):g\in G, x\in X\}^{\pm1}\subset E$
and let $Y$ be a set that bijects to this subset.  If $v=x_1x_2\ldots
x_n\in X^*$ then there is a $Y$-word $v'$ whose initial segments have
values $s(x_1), s(x_1x_2),\ldots, s(x_1\ldots x_n)$.  Let $L'$ be the
language
$$L'=\{v':v\in L\}\subset Y^*.$$

\th Proposition 2.2

The above language $L'$ is regular and has the following
fellow-traveller property: \Roster \ItemItem{(i)}There exists a
constant $K$ such that, if $w_1,w_2\in L'$ satisfy
$\pi(\eval{w_2})=\pi(\eval{y_1w_1y_2})$ with $y_1,y_2\in Y$ then
$y_1w_1y_2$ and $w_2$ $K$-fellow-travel in $E$.
\endRoster
Conversely, if $Y$ is a finite set which maps to a subset $\eval
Y=\eval Y^{-1}$ of $E$ and $L'\subset Y^*$ is a language with the
above property and:
\Roster
\ItemItem{(ii)} Evaluation maps $L'$ bijectively to the image of a section
$s\colon G\to E$,
\endRoster
\noindent then the projection $L$ of the language $L'$ to $G$ is a
biautomatic structure on $G$ and the cocycle defined by the section
$s$ is $L$-regular\endth

\pf Proof

We first show the fellow-traveller property for $L'$.  Let $Z$ be any
generating set for $A$.  Denote by $d_E$ and $d_G$ the word metrics in
$E$ and $G$ with respect to their generating sets $Y\cup Z$ and $X$.
It is readily established that $d_E(s(g),s(g'))=d_G(g,g')$ for all
$g,g'\in G$. Let $K$ be the fellow traveller constant for $L$.  Given
$w_1,w_2\in L'$ with $\pi(\eval{w_2})=\pi(\eval{y_1w_1y_2})$ and
$y_1,y_2\in Y$ then the fellow traveller property for $L$ tells us
that $d_G(\pi(\eval{y_1}w_1(t)),\pi(w_2(t)))\le K$ for all $t$. Thus
$d_E(s(\pi(\eval{y_1}w_1(t))),s(\pi(w_2(t))))\le K$.  But
$s(\pi(\eval{y_1}w_1(t)))$ differs from $s(\pi(\eval{y_1}))w_1(t)$ by
an element of $\sigma(X,G)$ and $s(\pi(\eval{y_1}))$ differs from
$\eval{y_1}$ by an element of $\sigma(G,X)$.  Also
$s(\pi(w_2(t)))=w_2(t)$.  Thus $d_E(\eval{y_1}w_1(t),w_2(t))\le
K+2K'$, where $K'$ is a bound on word-length in the sets $\sigma(X,G)$
and $\sigma(G,X)$.

We now show that the language $L'=\{w':w\in L\}$ is regular.  Denote
by $W$ the finite state automaton which has accepted language
$L$. Recall that $W$ may be regarded as a finite directed graph with
vertex set $S$, the elements of which are referred to as states.
There is a distinguished vertex, $\nu_0$, called the start state and a
distinguished subset of $S$, the elements of which are known as accept
states.  Each edge is labelled by an element of $X$, and each vertex
has exactly one outgoing edge for each element of $X$. The transition
function $\tau\colon S\times X\to S$ is given by setting
$\tau(\nu,x)=\nu'$ when there is an edge from $\nu$ to $\nu'$ labelled
by $x$.  A word in $X^*$ is accepted by $W$ precisely when it labels a
path beginning at the start state and ending at an accept state.  For
$x\in X$ and $a\in A$, let $W_{x,a}$ be the finite state automaton
which accepts the language $\{w\in L:\sigma(w,x)=a\}$.  Denote the
vertex set of $W_{x,a}$ by $S_{x,a}$ and the transition function for
$W_{x,a}$ by $\tau_{x,a}$. We form a finite state automaton for $L'$
by taking as vertex set the cartesian product $S\times(\prod S_{x,a})$
together with a single extra state $\varnothing$. The edge with
initial vertex $(\nu,\ldots,\nu_{x',a'},\ldots)$ labelled by
$s(x)a^{-1}$ has terminal vertex $\varnothing$, if $\nu_{x,a}$ is not
an accept state of $W_{x,a}$.  Otherwise, the terminal vertex is
$(\tau(\nu,x),\ldots,\tau_{x',a'}(\nu_{x',a'},x),\ldots)$. All edges
with initial vertex $\varnothing$ have terminal vertex $\varnothing$.
The start state is given by the vertex
$(\nu,\ldots,\nu_{x',a'},\ldots)$ which has $\nu=\nu_0$ and each
$\nu_{x',a'}$ the start state of $W_{x',a'}$.  The accept states are
given by vertices of the form $(\nu,\ldots,\nu_{x',a'},\ldots)$ with
$\nu\in S$, that is, any vertex which has an accept state as the first
coordinate. The finite state automaton we have defined has accepted
language $L'$.

For the converse statement suppose $L'$ is a language as in the
proposition.  Let $L$ be the projection of this language to a language
for $G$. Thus $L$ is the same formal language as $L'$ but with a
different evaluation map.  Then $L$ is certainly regular.  The
bisynchronous fellow-traveller property for $L$ is immediate from the
corresponding property (i) of $L'$. Thus $L$ is a biautomatic
structure.

Thus we only need to show that the cocycle $\sigma$ for the section
$s$ determined by $L'$ is regular.  The facts that $\sigma(Y,G)$ and
$\sigma(G,Y)$ are finite are easy consequences of the fellow-traveller
property (i) and we leave them to the reader.  For the rationality
statement note that the fellow-traveller property implies that the
language $\{(u,v)\in L'\times L': \eval u\eval
x\iota(-a+b)=\eval{v}\}$ is the language of a (synchronous) two-tape
automaton for any $x\in Y, a,b\in A$. Thus its projection onto its
first factor is regular.  Denote the image of $\eval x$ in $G$ by
$\hat x$. Then this projection is $\{u\in L':\exists v\in L',
\eval{ux}=\eval{v}\iota(a-b)\}$, the image of which in $G$ is $\{g\in
G:s(g)\eval x=s(g\hat x)\iota(a-b)\}$.  If we choose $b$ so $s(\hat
x)\eval x^{-1}=\iota(b)$ then this is $\{g\in G:s(g)s(\hat x)=s(g\hat
x)\iota(a)\}=\{g\in G:\sigma(g,x)=a\}$, so this set is rational.\qed

\th Corollary 2.3

If, in the situation of the above Proposition, $Z$ is a finite
$G$-invariant generating set for $A$ and we choose a $G$-invariant
biautomatic structure $L_A\subset Z^*$ on $A$ then $M=L'L_A$ is a
biautomatic structure on $E$.  (Structures $L_A$ as above always
exists --- cf.~\cite{ECHLPT} or \cite{NS2}.)
\endth

\pf Proof 

$M$ is certainly a regular language. Suppose $w_1,w_2\in L'$ and
$v_1,v_2\in L_A$ satisfy $\eval{xw_1v_1y}=\eval{w_2v_2}$ with $x,y\in
(Y\cup Z)$. Then $\eval xw_1$ $K$-fellow-travels $w_2$.  Hence
$d_E(\eval{v_1},\eval{v_2})\le K+1$, whence
$d_A(\eval{v_1},\eval{v_2})$ is bounded by some constant $c$ say.  It
follows that $v_1$ and $v_2$ are $cK_A$-fellow-travellers, where $K_A$
is the fellow-traveller constant for $L_A$. Now, since $L'$ is
injective it has a ``departure function'' (cf.~\cite{ECHLPT}), so
there exists a constant $\delta$ so that any subword $u$ of length at
least $\delta$ of an $L'$-word has $d(1,\eval u)>2K$.  Since $\eval
xw_1$ $K$-fellow-travels $w_2$, the lengths of $w_1$ and $w_2$ can
differ by at most $\delta$.  It follows easily that $\eval xw_1v_1$
fellow-travels $w_2v_2$ with constant $cK_A+c+1+K+\delta$.\qed

\th Corollary 2.4  

Let $A\to A'$ be an equivariant map of finitely generated abelian
groups with finite $G$-actions.  Suppose this map has finite kernel
and cokernel.  Let $L$ be a biautomatic structure on $G$. Then a class
in $H^2(G;A)$ is $L$-regular if and only if its image in $H^2(G;A')$
is $L$-regular.\endth

\pf Proof

The ``only if'' holds even if $A\to A'$ does not have finite kernel
and cokernel and is easy, so we shall just prove the ``if''.  A
homomorphism with finite kernel and cokernel is a composition of a
surjection with finite kernel and an injection with finite cokernel,
so it suffices to prove these two special cases.

Let $E$ and $E'$ be the virtually central
extensions determined by the cohomology classes in $H^2(G;A)$ and
$H^2(G;A')$ in question.  We have a commutative diagram
$$\matrix{
0&\longrightarrow&A&\longrightarrow&E
 &\longrightarrow&G&\longrightarrow&1\cr
&&\big\downarrow
&&\big\downarrow&&
\big\downarrow\rlap{$\vcenter{\hbox{$\scriptstyle =$}}$}\cr
0&\longrightarrow&A'&\longrightarrow&E' 
&\longrightarrow&G&\longrightarrow&1\cr}.
$$ 
Let $\sigma'$ be an $L$-regular cocycle representing the class in
$H^2(G;A')$.  Recall that $\sigma'$ is determined by some section
$s'\colon G\to E'$ and we have a regular language $L'\subset Y^*$ as
in Proposition 2.2 bijecting onto $s'(G)$, where $Y$ is some finite
set with an evaluation map to a symmetric subset of $E'$.

We first consider the case
that $A\to A'$ is surjective with finite kernel.  Then the same holds
for $E\to E'$.  Pick any lift of $Y\to E'$ to a map $Y\to E$ with
symmetric image and interpret $L'$ as a language on $E$.  Then $L'$
clearly satisfies the condition of Proposition 2.2, proving the
corollary in this case. 

Next suppose $A\to A'$ is injective with finite cokernel. Choose coset
representatives $a_1,\ldots,a_k\in A'$ for $A$ in $A'$.  Then
$\iota(a_1), \ldots,\iota(a_k)$ are coset representatives for $E$ in
$E'$. Let $c\colon E'\to \{a_1,\ldots,a_k\}$ be the map which picks
the coset representative.  Then the section $s\colon G\to E$ given by
$s(g)=s'(g)\iota(-c(s'(g))$ has cocycle
$\sigma(g,h)=\sigma'(g,h)+c(s'(gh))-c(s'(g))^{s'(h)}-c(s'(h))$.  This
is clearly weakly bounded and is easily seen to be regular.
\qed

Applying this corollary to the map $A\to A$ given by multiplying by a
non-zero integer shows:

\th Corollary 2.5

A cohomology class in $H^2(G;A)$ is ``{virtually $L$-regular}'' (that
is, some non-zero multiple can be represented by an $L$-regular
cocycle) if and only if it is regular.
\qed\endth

Now suppose $G$ is biautomatic with biautomatic structure $L$ and
$H<G$ is a subgroup of finite index.  Then there is an induced
biautomatic structure $L_H$ on $H$ which is unique up to equivalence.
Let $S$ be a set of right coset representatives for $H$ in $G$ and let
$r\colon G\to S$ be the map that takes an element to its coset
representative. The \Em{transfer map} $H^2(H;A)\to H^2(G;A)$ is
defined on the level of cocycles by the formula
$$T\sigma(g_1,g_2)=\sum_{y\in
S}\sigma(yg_1(r(yg_1))^{-1},r(yg_1)g_2(r(yg_1g_2))^{-1})^y.$$

\th Proposition 2.6

Suppose $H<G$ is a subgroup of finite index and $\sigma$ is an
$L_H$-regular cocycle on $H$ with coefficients in $A$.  Then
$T(\sigma)$ is an $L$-regular cocycle on $G$.\endth

\pf Proof

Since $T\sigma(g,x)=\sum_{y\in
S}\sigma(yg(r(yg))^{-1},r(yg)x(r(ygx))^{-1})^y$, the set $\{g\in
G:T\sigma(g,x)=a\}$ is the union over all sums of the form $\sum_{y\in
S}a_y=a$ of the sets $\bigcap_{y\in S}\{g\in G:
\sigma(yg(r(yg))^{-1},r(yg)x(r(ygx))^{-1})=a_y^{y^{-1}}\}$.  This is a
finite union of finite intersections, so it suffices to show that the
sets involved in the intersections are rational.  Now $\{g\in G:
\sigma(yg(r(yg))^{-1},r(yg)x(r(ygx))^{-1})=a_y^{y^{-1}}\}=\bigcup_{b\in
S} (\{g\in G:\sigma(ygb^{-1},bx(r(bx))^{-1}=a_y^{y^{-1}}\}\cap\{g\in
G:r(yg)=b\})$.  The set $\{g\in
G:\sigma(ygb^{-1},bx(r(bx))^{-1}=a_y^{y^{-1}}\}$ is a two-sided
translate of the rational set $\{g\in
G:\sigma(g,bx(r(bx))^{-1}=a_y^{y^{-1}}\}$ and is hence rational, while
$\{g\in G:r(yg)=b\}$ is a translate of a subgroup of finite index and
is hence rational.\qed

\th Corollary 2.7

If $H<G$ is of finite index then the restriction of a cohomology class
$x\in H^2(G;A)$ to $H^2(H;A)$ is $L_H$-regular if and only if $x$ is
$L$-regular.\endth

\pf Proof

The ``if'' is easy so we prove the ``only if.'' Thus, assume the
restriction of $x$ is regular.  Since the composition of restriction
and transfer $H^2(G;A)\to H^2(H;A)\to H^2(G;A)$ is multiplication by
the index $[G:H]$, it follows that the element $[G:H]x\in H^2(G;A)$ is
regular, so $x$ is virtually regular. Thus the result follows from
Corollary 2.5.\qed

\pf Proof of Theorem A

Corollary 2.3 is one direction of Theorem A in the
introduction.  To prove the other direction we appeal to the work of
Lee Mosher \cite{M}.  He proves that if a central extension $E$ of a
group $G$ has a biautomatic structure then so does $G$.  His main
argument is the construction of a language $L'$ satisfying the
conditions of Proposition 2.2 above, in the case of a central
extension
$$0\to\Bbb Z\to E\to G\to 1.$$  In particular, in this situation
Proposition 2.2 then says the cohomology class for the extension is
regular.

We first consider the case of a central extension
$$0\to A\to E\to G\to 1,$$ such that
$E$ has a biautomatic structure. Let $x\in H^2(G;A)$ be its cohomology
class.  Write $A$ as a direct sum of a finite group $F$ and copies of
$\bbZ$ as follows: $A=F\oplus\coprod_{i=1}^n \bbZ$. Then
$H^2(G;A)=H^2(G;F)\oplus\coprod_{i=1}^n H^2(G;\bbZ)$.  For each
$j=1,\ldots,n$ we can form $K_j=E/(F\oplus\coprod_{i\ne j}\bbZ)$ and
we have the induced extension
$$0\to\Bbb Z\to K_j\to G\to
1.\eqno{(*)}$$ Lee Mosher's results say firstly that $K_j$ is
biautomatic (since $E$ is a central extension of $K_j$) and therefore
secondly, via the above remarks, that the cohomology class of $(*)$ is
regular.  That is, the image of $x$ in the $j$-th summand
$H^2(G;\bbZ)$ of $H^2(G;A)$ is regular for each $j=1,\ldots,n$.  By
Corollary 2.4 the same is true for the image of $x$ in $H^2(G;F)$.  It
follows that $x$ is regular.

Now if the extension is only a virtually central extension we take $H$
to be the kernel of the action of $G$ on $A$ and consider the
restriction of our extension: $0\to A\to E_0\to H\to 1$. This is a central
extension, so we can apply the case just proven to it and then apply
Corollary 2.7 to complete the proof of Theorem A. \qed

\rk Remark 2.8

If one replaces ``biautomatic'' by ``automatic'' or ``asynchronously
automatic'' in the above discussion, then it is appropriate to replace
the concept of ``regular'' cocycle by a concept ``right regular''
obtained by dropping the condition that $\sigma(X,G)$ be finite.  The
analogs of the results 2.1--2.5 then go through, though we do not know
if the analog of Theorem A holds.
\endrk

\HH 3.~~Biautomatic structures for virtually central extensions of
Fuchsian groups

\th Theorem 3.1  

Let $G$ be a finitely generated Fuchsian group. Then any virtually
central extension of $G$ by a finitely generated abelian group has a
biautomatic structure.\endth

\pf Proof

We shall use the geodesic language $L$ with respect to any finite
generating set as a biautomatic structure on $G$. Let $A$ be any
finitely generated abelian group with finite $G$-action.  It suffices
to show that every class in $H^2(G;A)$ is $L$-regular.  By Corollary
2.7 we may replace $G$ by a subgroup of finite index as desired.  Thus
there is no loss of generality in assuming $G$ is torsion free and
acts trivially on $A$, so we will do so.  As in the previous section,
we can then split $A$ as the sum of copies of $\bbZ$ and a finite
group $F$. Any class in $H^2(G;F)$ is regular by Corollary 2.4, so it
suffices to prove that any class in $H^2(G;\bbZ)$ is regular. 

If $\bbH^2/G$ is non-compact then $G$ is free, so $H^2(G;\bbZ)=0$.
Thus assume that $\bbH^2/G$ is compact.  Then Gersten in \cite{G2} in
effect constructed a regular cocycle $\sigma$ representing the
generator of $H^2(G;\bbZ)$ (we give a different construction below).
Thus every element of $H^2(G;\bbZ)$ is regular.\qed

Gersten has informed us that his construction of the biautomatic
structure in the torsion free case will remain unpublished. We
therefore give a treatment of his result here for completeness.

Our construction is rather different from his and yields a regular
cocycle for a multiple of the generator of $H^2(G;\bbZ)$ when $G$ is a
closed surface group of genus $g>1$, rather than for the generator.

Fix a presentation $$G=\langle a_1,b_1,\ldots,a_g,b_g\mid
\prod_{i=1}^g[a_i,b_i]=1\rangle,$$ and let $P$ be the hyperbolic
$4g$-gon with angles $\pi/2g$ and sides labelled by the $a_i$ and
$b_i$ in such a way that the a word corresponding to a circuit of $P$
is the relator for the above presentation. Identifying corresponding
sides of $P$ gives a hyperbolic structure on the closed surface of
genus $g$, and there is a tessellation of $\bbH^2$ by copies of $P$
given by the universal cover of the surface. The 1-skeleton $\Gamma$
of this tessellation is the Cayley graph of $H$ with respect to the
generating set $X=\{a_1,b_1,\ldots,a_g,b_g\}^{\pm1}$. Suppose
$w=x_1\ldots x_n\in X^*$ is a word.  We consider a point moving along
the path of this word. The tangent vector at the point is well defined
except at vertices of the path.  As we pass from the $x_i$ edge to the
$x_{i+1}$ edge of the path the tangent vector swings through an angle
of $\theta_i$ with $-\pi<\theta_i\le\pi$ (here $\theta_i=\pi$ only
occurs if $w$ is non-reduced, namely $x_{i+1}=x_i^{-1}$).  Define an
integer $n(w)$ by $$n(w)=\sum_{i=1}^{n-1}{2g\over \pi}\theta_i.$$
Notice that $n(w^{-1})=-n(w)$ if (and only if) $w$ is a reduced word.
If $w$ is a closed path, and we set $\theta_n$ equal to the angle that
the tangent vector swings through from the $x_n$ edge to the
$x_1$-edge, then it is a standard result of hyperbolic geometry that
$$\sum_{i=1}^n\theta_i=A(w)+2\pi \tau(w),$$ where $A(w)$ is the
``signed area'' enclosed by $w$ and $\tau(w)\in\bbZ$ is the ``turning
number'' of $w$, that is, the total rotation number of the tangent
vector as it moves along the path (we measure this either by parallel
translating all the tangent vectors back to some fixed base point in
$\bbH^2$ or by following the motion of a point at infinity determined
by the moving tangent vector). Thus
$$\eqalignno{n(w)&=2gA(w)/\pi+4g\tau(w)-{2g\over\pi}\theta_n\cr
&=8g(g-1)N(w)+4g\tau(w) -{2g\over\pi}\theta_n,&(**)\cr}$$ where $N(w)$
is the signed number of copies of $P$ enclosed by $w$.

Let $$E_k=\langle A_1,B_1,\ldots,A_g,B_g, Z\mid Z\ \hbox{\rm central},
\prod_{i=1}^g[a_i,b_i]=Z^k\rangle,$$ where $k=8g(g-1)$. The central
extension $$0\to\bbZ\mathop{\to}\limits^\iota
E_k\mathop{\to}\limits^\pi G\to 1$$ where $\pi(A_i)=a_i$,
$\pi(B_i)=b_i$ and $\pi(Z)=1$, represents $k$ times a generator of
$H^2(G;\bbZ)$. We shall construct a section for which the
corresponding cocycle is regular.

Let $L\subset X^*$ be a language which bijects to $G$ and comprises
only geodesic words. For $w=x_1\ldots x_n\in L$ denote by $W=X_1\ldots
X_n$ the word obtained by replacing each $a_i^\pm$ by $A_i^\pm$ and
each $b_i^\pm$ by $B_i^\pm$. Define a section $s\colon G\to E_k$ by
$s(\eval{w})=\eval{W}Z^{-n(w)}$.  

\th Proposition 3.2

With the above definitions, the cocycle $\sigma$ defined by
$\iota(\sigma(g_1,g_2))=s(g_1)s(g_2)s(g_1g_2)^{-1}$ is a bounded
regular cocycle.
\endth

\pf Proof

Note that the number $N(w)$ in $(**)$ can also be described as
follows. Since $\eval w=1$, we can write $w$ in the free group on $X$
as $$w=\prod_{j=1}^ru_jr^{n_j}u_j^{-1},$$ where
$r=\prod_{i=1}^g[a_i,b_i]$.  Then $N(w)=\sum_{j=1}^rn_j$.  Now if $w$
is a word with $\eval w=1$ and $W$ is the corresponding word
in the $A_i$ and $B_i$ then equation $(**)$ implies that $$\eval
WZ^{-n(w)}=Z^{-4g\tau(w)+{2g\over\pi}\theta_n}.$$

We first show that the cocycle $\sigma$ is a bounded cocycle.  Let
$g_1,g_2,g_3\in G$ with $g_1g_2=g_3$ and let $w_1,w_2,w_3\in L$ be the
words representing them. Then $s(g_i)=\eval{W_i}Z^{-n(w_i)}$, so
$$\eqalign{\iota(\sigma(g_1,g_2))&= s(g_1)s(g_2)s(g_3)^{-1} \cr
&=\eval{W_1W_2W_3^{-1}}Z^{-n(w_1)-n(w_2)+n(w_3)}.}$$ Denote
$w=w_1w_2w_3^{-1}$.  It is not hard to see that the path determined by
$w$ has $|\tau(w)|\le 2$.  Denote by $\phi_1$ the angle between the
tangent vectors to $w$ at the last edge of $w_1$ and the first edge of
$w_2$. Similarly, $\phi_2$ denotes the angle from $w_2$ to $w_3$, and
$\phi_3$ the angle from $w_3$ to $w_1$. We choose these with
$-\pi<\phi_i\le\pi$. Since $w_3$ is reduced, we have
$n(w_3^{-1})=-n(w_3)$, so $n(w_1)+n(w_2)-n(w_3)$ differs from $n(w)$
just by ${2g\over\pi}(\phi_1+\phi_2)$.  Thus
$$\iota(\sigma(g_1,g_2))=
Z^{-4g\tau(w)+{2g\over\pi}(\phi_1+\phi_2+\phi_3)},$$
and it follows that $\sigma$ is a bounded cocycle.

To prove that the cocycle is regular we consider the above formula in
case $g_2=\eval x$, where $x$ is a generator.  The language
$\{(w_1,w_3)\in L\times L:\eval{w_3}=\eval{w_1x}\}$ is regular.
Suppose that $x_1\ldots x_m$ and $y_1\ldots y_n$ is a pair of words in
this language. The values of $\phi_1$ and $\phi_2$ are determined by
$x_m$ and $y_n$ respectively. Similarly, the value of $\phi_3$ is
given by $x_1$ and $y_1$.  It is not hard to
see that the turning number is also determined by the same data in
this case (namely $\tau(w)=-1$ if all three of the $\phi_i$ are
negative, and otherwise $\tau(w)=0$ if $\phi_3$ is negative or both
$\phi_1$ and $\phi_2$ are negative, and $\tau(w)=1$ in all other
cases). It is clear that these data can be checked by finite state
automata, and therefore one can construct a finite state automaton
which will accept the language $\{w\in
L:\sigma(\eval{w},\eval{x})=a\}$.  \qed

\HH 4.~~Questions

S. Gersten, in \cite{G1}, shows that if a central extension $E$ of a
bicombable group $G$ by a finitely generated abelian group $A$ is
given by a bounded cocycle then $E$ is bicombable.  His argument is
the same as the argument of our section 2 --- the only difference
being that regularity of languages is not important.  It follows that
his result is valid even if the cocycle is only weakly bounded.  A
natural question therefore is whether a weakly bounded cohomology
class on a finitely generated group is always bounded.  

This question is also relevant to quasi-isometry. Gersten shows that
if the cocycle is bounded then $G\times A$ is quasi-isometric to $E$,
but it is again not hard to see that weakly bounded suffices.  In fact
in this case only the condition that $\sigma(G,X)$ is bounded is
needed --- the map $g\times a\mapsto s(g)\iota(a)$ then gives a
quasi-isometry.  Moreover, if a quasi-isometry $G\times A$ to $E$
exists such that the composition $G\times\{1\}\to E\to G$ is a
quasi-isometry then the central extension is determined by a cocyle
with $\sigma(G,X)$ bounded.  But we know no example of a cohomology
class for a group which is represented by such a cocycle and is not
bounded.

  Thurston has claimed (unpublished) that central extensions of
word-hyperbolic groups by finitely generated abelian groups are
automatic.  Are they in fact biautomatic?  In fact, might every
2-dimensional cohomology class on a word-hyperbolic group be
representable by a bounded regular 2-cocycle?

\Bib        Bibliography

\rf {ECHLPT} D.~B.~A. Epstein, J.~W. Cannon, D.~F. Holt, S.~V.~F. Levy,
M.~S. Patterson, and  W.~P. Thurston, Word processing in groups, Jones and
Bartlett, 1992.

\rf {G1} S.~M.~Gersten, Bounded cocycles and combings of groups, 
Int. J. of Algebra and  Computation {\bf 2} (1992), no. 3, 307--326.

\rf {G2} S.~M.~Gersten, Bounded cohomology and combings of groups,
Preprint, Version 4.1.

\rf {M} L.~Mosher, Central quotients of biautomatic groups, Preprint.

\rf {NS1} W.~D.~Neumann and M.~Shapiro, Equivalent automatic structures and
their boundaries, Int. J. of Algebra and  Computation {\bf 2} (1992), no. 4,
443--469.

\rf {NS2} W.~D.~Neumann and M.~Shapiro, Automatic structures, rational growth,
and geometrically finite groups. Invent. Math. (to appear).

\endBib

\Coordinates
Department of Mathematics
The University of Melbourne
Parkville, VIC 3052, Australia

neumann@maths.mu.oz.au
ldr@maths.mu.oz.au
\endCoordinates
\bye